\documentclass[draft]{Dekker2}
\usepackage{latexsym}
\usepackage{amsfonts}
\usepackage{times}
\newcommand{\nab}{{\vec{\nabla}}}
\newcommand{\xbo}{{\vec{x}}}
\newcommand{\ybo}{{\vec{y}}}
\newcommand{\vbo}{{\vec{v}}}
\newcommand{\ubo}{{\vec{u}}}
\newcommand{\taubo}{{\vec{\tau}}}
\newcommand{\xnk}{{\vec{X}_N^k}}
\newcommand{\xnl}{{\vec{X}_N^l}}
\newcommand{\xnj}{{\vec{X}_N^j}}
\newcommand{\bk}{{\vec{B}^k}}
\newcommand{\snr}{{S_{N,r}}}
\newcommand{\hnr}{{h_{N,r}}}
\newcommand{\dnr}{{d_{N,r}}}
\newcommand{\dnq}{{d_{N,q}}}
\newcommand{\snq}{{S_{N,q}}}
\newcommand{\harq}{{\hat{a}_{rq}}}
\newcommand{\haql}{{\hat{a}_{ql}}}
\newcommand{\eqlr}{{e_{qlr}}}
\newcommand{\erql}{{e_{rql}}}
\newcommand{\heqlr}{{\hat{e}_{qlr}}}
\newcommand{\herql}{{\hat{e}_{rql}}}
\newcommand{\heqrl}{{\hat{e}_{qrl}}}
\newcommand{\anrq}{{a_{N,rq}}}
\newcommand{\anrqkj}{{a_{N,rq}^{kj}}}
\newcommand{\anql}{{a_{N,ql}}}
\newcommand{\anrqs}{{a_{N,rq}^{*,k}}}
\newcommand{\anrqskj}{{a_{N,rq}^{*,kj}}}
\newcommand{\anqrsjk}{{a_{N,qr}^{*,jk}}}
\newcommand{\anqls}{{a_{N,ql}^{*,k}}}
\newcommand{\hbeta}{{\hat{\beta}}}
\newcommand{\hw}{{\hat{W}}}
\newcommand{\hm}{{\hat{M}}}

\newcommand{\ta}{{a}}
\newcommand{\han}{{\hat{\alpha}_N}}
\newcommand{\rel}{{\mathbb{R}}}
\newcommand{\nat}{{\mathbb{N}}}
\newcommand{\bea}{\begin{eqnarray}}
\newcommand{\eea}{\end{eqnarray}}
\newcommand{\bef}{\begin{figure}}
\newcommand{\enf}{\end{figure}}
\newcommand{\ball}{\begin{array}{ll}}
\newcommand{\bal}{\begin{array}{l}}
\newcommand{\ea}{\end{array}}
\newcommand{\bde}{\begin{description}}
\newcommand{\ede}{\end{description}}

\newcommand{\sur}{{\sum_{r=1}^R}}
\newcommand{\sul}{{\sum_{l=1}^R}}
\newcommand{\suq}{{\sum_{q=1}^R}}
\newcommand{\sukrt}{{\sum_{k\in M(N,r,t)}}}
\newcommand{\sujrt}{{\sum_{j\in M(N,r,t)}}}

\newcommand{\sukrs}{{\sum_{k\in M(N,r,s)}}}
\newcommand{\sukqs}{{\sum_{k\in M(N,q,s)}}}
\newcommand{\limn}{{\lim_{N\rightarrow\infty}}}
\newcommand{\halb}{{1\over 2}}
\newcommand{\oon}{{1\over N}}

\newcommand{\bml}{\begin{list}{\textbf{(M\arabic{model})}}{\usecounter{model}\setlength{\leftmargin}{9mm}}}
\newcommand{\bcl}{\begin{list}{\textbf{(C\arabic{model})}}{\usecounter{model}\setlength{\leftmargin}{9mm}}}
\newcounter{model}

\begin{document}

\setcounter{page}{1}
\jname{Stoch.\ An.\ Appl.}
\jvol{}
\jissue{}
\jyear{2003}

\webslug{www.dekker.com}
\cpright{Marcel Dekker, Inc.}

\title[A rigorous derivation of Smoluchowski's equation]{A rigorous derivation of Smoluchowski's equation in the moderate limit}

\author[Gro\ss kinsky et al.]{S.\ Gro\ss kinsky$^1$, C.\ Klingenberg$^2$ and K.\ Oelschl\"ager$^3$\affiliation{$^1$Zentrum Mathematik M5, Technische Universit\"at M\"unchen, 85747 Garching bei M\"unchen, Germany\\
$^2$Institut f\"ur Angewandte Mathematik, Universit\"at W\"urzburg, Am Hubland, 97074 W\"urzburg, Germany\\
$^3$Institut f\"ur Angewandte Mathematik, Universit\"at Heidelberg, Im Neuenheimer Feld 294, 69120 Heidelberg, Germany}}

\abstract{Smoluchowski's equation is a macroscopic description of a many particle system with coagulation and shattering interactions. We give a microscopic model of the system from which we derive this equation rigorously. Provided the existence of a unique and sufficiently regular solution of Smoluchowski's equation, we prove the law of large numbers for the empirical processes. In contrast to previous derivations we assume a moderate scaling of the particle interaction, enabling us to estimate the critical fluctuation terms by using martingale inequalities. This approach can be justified in the regime of high temperatures and particle densities, which is of special interest in astrophysical studies and where previous derivations do not apply.\\}

\makechaptertitle

\section{Introduction}
We consider a system of dust particles of $R\in\nat$ different masses $m_1 ,\ldots ,m_R$, embeded in a $d$-dimensional hot gas. Particles of size $r\in\{ 1,\ldots ,R\}$ are drifting according to the velocity field $\vbo_r :\rel^d \times\rel_0^+ \rightarrow\rel^d$ with a superimposed Brownian motion with diffusion constant $\sigma_r \in\rel^+$. Two particles of size $r$ and $q$ collide with rate $\harq :\rel^d \times\rel_0^+ \rightarrow\rel_0^+$. The material coefficients $\herql\in\nat_0$ determine the number of particles of size $l=1,\ldots ,R$ produced by that collision event, deciding for coagulation or shattering events. In this first model we take a macroscopic viewpoint, where two colliding particles occupy the same position in space-time $(\xbo,t)$, where the function $\harq$ is evaluated. A complete description of the above model is given by Smoluchowski's equation \cite{smol}, in our case a system of reaction diffusion equations for the particle densities $s_r :\rel^d \times\rel_0^+ \rightarrow\rel_0^+$, $r=1,\ldots ,R$ (cf.\ \cite{yorke1}, sect.\ 2) with initial conditions $s_r^0 :\rel^d \rightarrow\rel_0^+$:
\bea
\quad\partial_t s_r (\xbo ,t)&=&-\nab\cdot\big(\vbo_r (\xbo ,t)\, s_r (\xbo ,t)\big) +{1\over 2}\sigma_r^2 \,\Delta s_r (\xbo ,t)\nonumber\\
& &-s_r (\xbo ,t)\suq\harq (\xbo ,t)\, s_q (\xbo ,t)+ \halb \sum_{q,l=1}^R \haql (\xbo ,t)\, \heqlr\, s_q (\xbo ,t)\, s_l (\xbo ,t)\nonumber\\
s_r (\xbo ,0)&=&s_r^0 (\xbo) \quad\mbox{for all }r=1,\ldots ,R\hspace{66mm}\ \mbox{(MA)}\nonumber
\eea
We suppose that for any size $r$ the particles consist of several atoms of size $1$ and we set $m_1 =1$. Moreover, the masses of the particles are ordered as $m_1 < m_2 <\ldots <m_R$. The conservation of the total mass $M(t)=\sur\int_{\rel^d} m_r \, s_r (\xbo,t)\, d^d x$ of the system under the above dynamics, i.e.\ $M(t)=M(0), t\geq 0$, is assured by
\bea\label{masscons}
\sul m_l \,\herql = m_r +m_q\qquad\mbox{and}\qquad\harq (\xbo,t)= \hat{a}_{qr} (\xbo,t)\ ,
\eea
for all $\xbo\in\rel^d ,\ t\in\rel_0^+ ,\ r,q\in\{ 1,\ldots ,R\}$. We also assume that $\herql$ is symmetric in $r$ and $q$, i.e.\ $\herql =\heqrl$.\\
The above model is commonly used to describe reaction diffusion systems and there have been rigorous approaches to identify equation (MA) as the limit dynamics of a suitable many particle system. These derivations are restricted to one-dimensional systems \cite{arnold}, the spatially homogeneous case \cite{norris} or a spatially discretized microscopic model \cite{demasi,guias}. In \cite{lang} there is a derivation accounting for the full space dependence of the problem, using the Boltzmann-Grad limit which is applicable for very small particle densities. In \cite{guias,norris} existence and uniqueness of a solution of Smoluchowski's equation are also studied.\\
In this paper we give a microscopic particle model (MI) in section II, from which we rigorously derive (MA) in the spatially inhomogeneous (general) case without space discretization or restrictions on space dimension. Our many particle system properly describes an astrophysical system recently studied in \cite{jones,yorke1}, which is explained in section V.A. It corresponds to a situation of high gas temperatures and particle densities, which is not covered by the derivation in the Boltzmann-Grad limit \cite{lang}.\\
In this regime the dominating particle interactions are shattering collisions, so it is justified to neglect coagulation events. That means that the mass of each of the two interaction partners may not increase by the collision, but they are shattered into fragments of smaller or equal mass. This constitutes a constraint on the material coefficients $\herql$ given in (\ref{micrcons}), which is important to ensure compatibility with the microscopic particle model. Our main theorem in section III states the convergence of the empirical processes (\ref{snrdef}) to a solution of (MA) and is proved in section IV. Before giving a short conclusion in the last section we also discuss two apparent generalizations of the microscopic model (MI).\\
The most important feature of our approach is the moderate scaling of the collision interaction, which is introduced in section II.B (M3) and discussed on a physical level in section V.A. It enables us to use a technique developed by K.\ Oelschl\"ager \cite{oel1}, which was previously applied to derive the porous medium equation \cite{oel4}, or in the description of aggregation phenomena in biological populations \cite{morale1}, \cite{stevens1}. With this technique we are able to derive Smoluchowski's equation in the spatially inhomogeneous form (MA), in a regime where the previous approaches cannot be applied.

\section{Microscopic particle model}
Given the macroscopic model of section I we present a corresponding microscopic many particle system. The most important modeling assumptions are marked by (M1) to (M4) and are discussed in sections V.A and V.B.

  \subsection{Dynamics without interaction}
Let $N_r (t)$ be the number of particles of species $r\in \{1,\ldots ,R\}$ and $N(t) =\sur N_r (t)$ the number of all particles at time $t\in\rel_0^+$. The system size $N$ is characterized by the number of atoms of mass $m_1 =1$ at time $t=0$:
\bea\label{ndef}
N=\left.\left(\sur N_r (0)\, m_r \right)\!\right/\! m_1 =\sur N_r (0)\, m_r
\eea
Let $M(N,t)\subset\nat$ be the set of all particles and $M(N,r,t)\subset M(N,t)$, $r=1,\ldots ,R$, the subsets of particles of species $r$ at time $t$, where each particle is identified with a unique integer number.
\bml
\item The particles are considered to be point masses with positions $\xnk (t)\in\rel^d$, $k\in M(N,t)$, at time $t$ in a system of size $N$. Each particle of species $r\in\{ 1,\ldots ,R\}$ is given the rescaled mass $m_{N,r} =m_r /N$, which keeps the initial total mass $M_N (0)=\sur m_{N,r} N_r (0)=1$ independent of the system size according to (\ref{ndef}).
\item Neglecting the hydrodynamic drag interaction between gas and particles, we consider the latter to move according to the given velocity fields $\vbo_r$ and Brownian motion with diffusion constants $\sigma_r$, $r=1,\ldots ,R$, introduced in the macroscopic equation (MA). 
\end{list}
Between two subsequent collision events the system at time $t$ is then described by $N(t)$ uncoupled stochastic differential equations:
\bea\label{micr1}
d\xnk (t)=\vbo_r (\xnk(t),t)\, dt+\sigma_r d\bk (t)\ ,\quad k\in M(N,r,t)\ ,\quad r=1,\ldots ,R
\eea
The $(\bk (t))_{t\in\rel_0^+}$, $k\in\nat$, are independent Wiener processes modelling the Brownian motion of the particles. We always assume the existence of a filtration $(\mathcal{F}_t )_{t\geq 0}$, with respect to which the stochastic processes under consideration are adapted (cf.\ \cite{oel1}, sect.\ 2.B) and which fulfills the usual conditions \cite{karatzas}.\\
The particle interaction is described by suitable changes of the sets $M(N,r,t)$ and is explained in the next subsection. A microscopic quantity comparable to the particle density in (MA) is given by the measure-valued, empirical processes:
\bea\label{snrdef}
\snr :\rel_0^+ \rightarrow\mathcal{M}(\rel^d )\ ,\quad\snr (t):=\oon\sukrt\delta_{\xnk(t)} \ ,\qquad r=1,\ldots ,R,
\eea
where $\mathcal{M}(\rel^d )$ denotes the space of positive, finite measures on $\rel^d$ and $\delta_\xbo$ is the Dirac measure concentrated in $\xbo\in\rel^d$. $\snr$ describes the time-evolution of the spatial distribution of particles within the subpopulation of species $r$. It is known by the law of large numbers that the empirical distribution of $N$ independent, identically distributed random variables converges to their probability distribution in the limit $N\rightarrow\infty$. In this paper we prove the convergence for stochastic processes which are not independent for times $t>0$, due to the particle interaction.

  \subsection{Description of the particle interaction}
Due to (M1) we have to specify a model for the `collision' interaction of two point particles.
\bml
\setcounter{model}{2}
\item We take a stochastic model determined by a rate depending on the distance of the interaction partners $k$ and $l$. The scaling of this rate is given by
\bea\label{rate}
\oon\, W_N (\xnk(t)-\xnl(t))\ ,\quad\mbox{where }W_N (\xbo)=\alpha_N^d W_1 (\alpha_N \xbo)\ ,
\eea
with $\alpha_N =N^{\beta /d}$ and a moderate scaling parameter $0<\beta <1$. We assume that $W_1$ is symmetric and positive with $\| W_1\|_1 =1$. It follows that $\| W_N\|_1 =1$ for all $N\in\nat$ and\linebreak $\limn W_N (.-\xbo)=\delta_\xbo$ for all $\xbo\in\rel^d$ in the sense of distributions.
\end{list}
In contrast to the usual hydrodynamic scaling with $\beta =1$ this leads to a microscopically large interaction volume. This assumption is motivated and justified in a physical context in section V.A.
\bml
\setcounter{model}{3}
\item Instead of considering pair interactions (see sect.\ V.B) we assume that every particle $k\in M(N,r,t)$ interacts with an effective field of all other particles of species $q$ with rate
\bea\label{adef}
\!\!\!\!\anrq (\xnk (t),t):=\min\Big\{ C_a \!\!&,&\!\harq (\xnk(t),t)\nonumber\\
& &\!\Big((\snq(t) *W_N )(\xnk(t))\! -\!\delta_{r,q} W_N (\vec{0})/N \Big)\Big\} ,
\eea
where $\harq$ is the macroscopic collision rate given in section I.
\end{list}
In (\ref{adef}) we used the generalized convolution product
\bea\label{convdef}
\qquad\snq (t)*W_N :=\int_{\rel^d}\! W_N (.-\xbo)\, \snq (t)(d^d x)=\oon\!\sum_{l\in M(N,q,t)}\!\!\! W_N (.-\xnl(t)).
\eea
By substraction of the term including Kronecker's delta in (\ref{adef}) self-interaction is excluded. The rate is bounded uniformly in $N$ by a suitable constant $C_a$, which is specified in condition (C5) in section III.B. This cut-off prevents diverging interaction rates due to high particle concentrations in the limit $N\rightarrow\infty$. Each possible interaction event is described by a jump process
\bea\label{jump}
\anrqs (t):=\beta_{N,rq}^k \left(\int_0^t \chi_{M(N,r,s)} (k)\, \anrq (\xnk(s),s)\, ds\right)\quad\in\quad\{ 0,1\} ,
\eea
where $\beta_{N,rq}^k :\rel_0^+ \rightarrow\nat_0$ are independent standard Poisson processes with a transformed time argument in the brackets (cf. \cite{oel1}) and $\chi_A \in\{ 0,1\}$ is the indicator function of the set $A$.\\
The process $\anrqs$ jumps from $0$ to $1$ at some time $t\geq 0$ if particle $k$ exists in $t-$, belongs to species $r$ and interacts with a particle of species $q$ at time $t$. After the interaction the number $k$ is removed from the sets $M(N,r,t)$ and $M(N,t)$. The mass of particle $k$ is distributed on the interaction products according to the microscopic material coefficient $\erql\in\nat_0$. The latter fulfills conservation of mass and is related to its macroscopic counterpart $\herql$ in the following way:
\bea\label{micrcons}
\sul m_l \erql =m_r \ ,\quad\herql =\heqrl =\erql +e_{qrl} \quad\mbox{for all }r,q,l=1,\ldots ,R.
\eea
We note that this also constitutes a condition on $\herql$, corresponding to the absence of coagulation mentioned in section I. The particles resulting from the interaction are located at $\xnk(t)$ and obtain new numbers starting with $\max\{ p\in M(N,s):s\leq t\} +1$, which were previously not assigned to any particle. These numbers are added to $M(N,t)$ and the subsets corresponding to the various species. We note that any process $\anrqs$ only jumps once, since after that jump the respective particle $k$ disappears, i.e.\ $\chi_{M(N,r,t+.)}(k)\equiv 0$.

  \subsection{Complete description of the model}
Using a generalized $L^2$-scalar product we can formulate the time evolution of the empirical processes in a weak sense. For all $f\in C^2_b (\rel^d ,\rel )$ and $r=1,\ldots ,R$ we have
\bea\label{micr2}
\langle S_{N,r}(t),f\rangle &:=&\int_{\rel^d} f(\xbo)\, \snr(t)(d^d x)={1\over N}\sukrt f(\xnk (t))=\nonumber\\
&=&{1\over N}\sum_{k\in M(N,r,0)} f(\xnk (0))+\int_0^t d\left(\oon\sukrs f(\xnk (s))\right) .
\eea
Inserting the expression for $d\xnk (s)$ from equation (\ref{micr1}) and using It\^o's formula \cite{karatzas} we get:
\bea\label{micr4}
\langle S_{N,r}(t),f\rangle&=&{1\over N}\sum_{k\in M(N,r,0)} f(\xnk (0))+{\sigma_r \over N}\int_0^t \sukrs \nab f(\xnk(s))\cdot d\bk (s)\nonumber\\
& &+{1\over N}\int_0^t \sukrs\left( \nab f(\xnk(s))\cdot\vbo_r (\xnk(s),s)+{\sigma_r^2 \over 2}\Delta f(\xnk(s))\right) ds\nonumber\\
& &-{1\over N}\suq\int_0^t \sukrs f(\xnk(s))\, \anrqs (ds)\nonumber\\
& &+{1\over N}\sum_{q,l=1}^R \int_0^t \sukqs f(\xnk(s))\,\eqlr\,\anqls (ds)
\eea
The first integral term describes the stochastic fluctuations of the particle positions and the second one particle transport and diffusion, resulting from the interaction free description (\ref{micr1}). The next two terms consider the change of the sets $M(N,r,s)$ in (\ref{micr2}) due to the loss of particles of species $r$ after interactions with others, and the gain of such particles from products of other interactions. We separate the fluctuation terms due to stochasticity in the free particle dynamics and the interaction in stochastic integrals. So we get for all $f\in C^2_b (\rel^d ,\rel )$ and $r=1,\ldots ,R$ the complete description of our microscopic model:
\bea
\langle S_{N,r}(t),f\rangle &=&\langle S_{N,r}(0),f\rangle +\int_0^t \left\langle S_{N,r}(s),\nab f\cdot\vbo_r (.,s)+ {\sigma_r^2\over 2} \Delta f\right\rangle ds\nonumber\\
& &-\suq\int_0^t \!\langle S_{N,r}(s),f\, \anrq (.,s)\rangle\, ds+\!\!\sum_{q,l=1}^R \int_0^t \!\langle S_{N,q}(s),f\,\eqlr\,\anql (.,s)\rangle\, ds\nonumber\\
& &+M_{N,r}^1 (f,t)+M_{N,r}^{2a} (f,t)+M_{N,r}^{2b} (f,t)\ ,\quad\mbox{with}\nonumber\\
& &\nonumber\\
M_{N,r}^1 (f,t)&=&{\sigma_r \over N}\int_0^t \sukrs \nab f(\xnk(s))\cdot d\bk (s)\nonumber\\
M_{N,r}^{2a} (f,t)&=&-{1\over N}\suq\int_0^t \sukrs f(\xnk(s))\Big( \anrqs (ds)-\anrq (\xnk(s),s)ds\Big)\nonumber\\
M_{N,r}^{2b} (f,t)&=&{1\over N}\sum_{q,l=1}^R \int_0^t \sukqs f(\xnk(s))\,\eqlr\Big( \anqls (ds)-\anql (\xnk(s),s)ds \Big)\nonumber
\eea
and initial conditions $\snr (0)=N^{-1} \sum_{k\in M(N,r,0)}\delta_{\xnk(0)}$.\hfill (MI)\\
\\
This set of equations combines all features mentioned in the preceding two subsections and is used to derive the macroscopic model (MA), shown in the next section.

\section{Derivation of Smoluchowski's equation}
We show how to obtain (MA) heuristically from our microscopic particle model (MI), leading us to a proper formulation of the main theorem.

  \subsection{Heuristic derivation of the macroscopic equation}
The empirical processes $\snr$ are defined as solutions of (MI). For this subsection we assume that for every $r=1,\ldots ,R$ they converge to limit processes $S_r :[0,T] \rightarrow\mathcal{M}(\rel^d )$ on a compact time interval $[0,T]$ in a yet unspecified sense. The limit processes are assumed to be absolutely continuous with respect to Lebesgue measure on $\rel^d$ and therefore have densities $s_r$, which should be in $C_b^2 (\rel^d \times\rel_0^+ ,\rel_0^+ )$. With the generalized scalar product defined in (\ref{micr2}) we therefore have $\langle S_r (t),g(.,t)\rangle =\langle s_r (.,t),g(.,t)\rangle$ for all $g\in C_b (\rel^d \times [0,T],\rel )$. We also assume the validity of conditions (C1) to (C7) given in the next subsection.\\
In section IV.D we get the following for the stochastic integrals in (MI) for any $T>0$:
\bea\label{erg1}
\limn E\left[\sup_{t\in [0,T]} |M_{N,r}^i (f,t)|\right] =0\quad\mbox{for }i=1,2a,2b\ ,
\eea
so the fluctuation terms asymptotically vanish in any compact time interval and the limit equation is supposed to be deterministic (see (\ref{215})). The convergence of the $\snr$ should be sufficiently strong to assure the following:
\bea\label{limcond}
\limn E\left[\langle\snr(t),g(.,t)\rangle\right] &=&\langle s_r (.,t),g(.,t)\rangle\nonumber\\
\limn E\left[\langle\snr(t),(\snq(t)*W_N )\, g(.,t)\rangle\right] &=&\langle s_r (.,t),s_q (.,t)\, g(.,t)\rangle ,
\eea
for all $r,q=1,\ldots ,R$, $t\in[0,T]$ and $g\in C_b (\rel^d \times [0,T],\rel^d )$. The first condition assures the convergence of the drift and diffusion term in (MI) and the second one is needed for the interaction terms. We formally substitute the above limits into (MI) and notice that the self interaction term in (\ref{adef}) vanishes for $N\rightarrow\infty$. Therefore we get the following deterministic integral equation for all test functions $f\in C_b^2 (\rel^d ,\rel)$, $t\in [0,T]$ and $r=1,\ldots ,R$:
\bea\label{215}
\langle s_r (.,t),f\rangle &=&\langle s_r (.,0),f\rangle +\int_0^t ds\left\langle s_r (.,s),\nab f\cdot\vbo_r (.,s)+ {1\over 2}\sigma_r^2 \Delta f\right\rangle\nonumber\\
& &-\suq\int_0^t ds\big\langle s_r (.,s),f\, \harq (.,s)\, s_q (.,s)\big\rangle\nonumber\\
& &+\sum_{q,l=1}^R \int_0^t ds\,\big\langle s_q (.,s),f\, \eqlr\,\haql (.,s)\, s_l (.,s)\big\rangle
\eea
After partial integration in the transport and diffusion terms one immediately recognizes this as a weak version of Smoluchowski's equation. Using (\ref{masscons}) and (\ref{micrcons}) it is easy to get the last line in the form (MA).\\
Therefore we showed that, assuming the empirical processes converge, their limit densities fulfill a weak form of Smoluchowski's equation. In the next subsection we explain how to prove this convergence in an appropriate rigorous limit sense, which can be seen from (\ref{limcond}) to be of $L^2$-type.

  \subsection{Convergence theorem}
To formulate the convergence theorem we use the following distance function between the empirical processes (MI) and the solution of Smoluchowski's equation (MA) specified in (C3) below:
\bea\label{dnrdef}
\dnr(.,t):=h_{N,r} (.,t)-s_r (.,t)\ ,\quad\mbox{where}\quad h_{N,r} (.,t):=\snr (t)*\hw_N ,
\eea
for all $r=1,\ldots ,R$ and $t\in\rel_0^+$. The convolution kernel $\hw_N$ smooths out the empirical processes and obeys the following regularity conditions:
\bcl
\item $\hw_N$ is a different scaling of the interaction function $W_N$ and both have to fulfill:
\bea
\hw_N (\xbo)=\han^d W_1 (\han\xbo)\qquad\quad\ &\mbox{and}&\quad W_N (\xbo) =\alpha_N^d W_1 (\alpha_N \xbo)\ ,\quad\mbox{where}\nonumber\\
\han =N^{\hbeta /d},\ 0<\hbeta<{d\over d+2}\quad&\mbox{and}&\quad\alpha_N =N^{\beta /d}\mbox,\ 0<\beta <{\hbeta\over d+1} .\nonumber
\eea
\end{list}
The scaling parameter $\hbeta$ plays no role in the dynamics of the many-particle system. However, by the above assumptions some restrictions on the parameter $\beta$ determining the moderate interaction are introduced.
\bcl
\setcounter{model}{1}
\item The unscaled function $W_1 \in L^1 \cap C_b^2 (\rel^d ,\rel)$ is symmetric, positive and standardized, i.e.\ $\| W_1 \|_1 =1$. We also need $\int_{\rel^d} |\xbo |W_1 (\xbo) d^d x <\infty$ and the Fourier transform $\tilde{W}_1$ has to fulfill:
\bea
&a)&\tilde{W}_1 \in C_b^2 (\rel^d)\nonumber\\
&b)&|\tilde{W}_1 (\taubo )|\leq C\exp (-C'|\taubo |)\nonumber\\
&c)&|\Delta \tilde{W}_1 (\taubo)|\leq C(1+|\taubo |^2 )|\tilde{W}_1 (\taubo )|\nonumber\\
&d)&v\mapsto |\tilde{W}_1 (v\taubo )|,\ v\geq 0,\mbox{ monotonicly decreasing for all fixed }\taubo\in\rel^d\nonumber
\eea
\end{list}
A Gaussian probability density is an example for $W_1$ which obeys these conditions. To the knowledge of the authors there is no proof of the existence of a sufficiently smooth solution of the macroscopic equations, therefore we have to assume the following:
\bcl
\setcounter{model}{2}
\item There exists a positive, unique $C_b^2$ -solution $(s_1 ,\ldots ,s_R)$ of Smoluchowski's equation (MA) in the time interval $[0,T^* ]$ for some positive $T^*$. The functions $s_r (.,t)$ and their partial derivatives are $L^2 (\rel^d ,\rel_0^+ )$-bounded uniformly in $t\in [0,T^* ]$.
\item The macroscopic collision rate $\harq(.,t)$ given in section I should be Lipschitz continuous, bounded and fulfill the conditions (\ref{masscons}) for all $t\in [0,T^* ]$. The macroscopic material coefficient $\herql$ should obey condition (\ref{masscons}) and together with its microscopic counterpart $\erql$ given in section II.B, it should fulfill (\ref{micrcons}) and be symmetric in $r$ and $q$.
\item The upper bound $C_a$ for the microscopic interaction rates (\ref{adef}) is given so that the limit equation is not affected, $C_a >\max_{r,q\in\{ 1,\ldots ,R\} } \sup_{t\in [0,T^* ]} \|\harq (.,t)s_q (.,t)\|_\infty$.
\item The velocity fields of the different particle species have to fulfill
\bea
\vbo_r \in C_b^1 (\rel^d \times [0,T^*] ,\rel^d )\quad\mbox{for all }r=1,\ldots ,R.\nonumber
\eea
\item The diffusion constants of all particle species have to be positive, i.e.\ $\sigma_r >0$ for all $r=1,\ldots ,R$.
\end{list}
We note that our proof only applies if all particles are Brownian. Now we are ready to formulate our main convergence result.\\[2mm]
\textbf{Theorem.}\it\ \ With conditions (C1) to (C7) and $\limn E\left[\sur\|\dnr(.,0)\|_2^2 \right] =0$ it is\rm
\bea\label{stat1}
\limn E\left[\sur \sup_{t\in [0,T^* ]} \|\dnr(.,t)\|_2^2 \right] =0.\\
\nonumber
\eea
Convergence at time $t=0$ is given if the initial conditions of (MA) and (MI) are compatible. One possibility is to take the particle positions $\xnk(0)$, $k\in M(N,r,0)$ as independent, identically distributed random variables with suitably normalized densities $s_r^0 /\sur\langle s_r^0 ,1\rangle$ for all $r=1,\ldots ,R$. For discussion of this point see \cite{oel4} (sect.\ 4B).\\
To formulate the result without the smoothing convolution kernel $\hw_N$ we introduce a metric on $\mathcal{M} (\rel^d )$ by
\bea\label{metrdef}
D(\mu ,\nu )&:=&\sup \{ |\langle\mu -\nu ,f\rangle |:f\in \Omega_D\}\quad\mbox{for all }\mu ,\nu\in\mathcal{M}(\rel^d )\mbox{ and}\nonumber\\
\Omega_D&:=&\{ f\in C_b^1 \cap L^2 (\rel^d ,\rel)\, :\,\|f\|_\infty +\|\nab f\|_\infty +\|f\|_2 \leq 1\} .
\eea
This quantifies a distance between the empirical processes defined in (MI) and the processes $S_r (t):=\int_. s_r (\xbo ,t)d^d x$ given by the solution (C3) of the macroscopic equation. As the theorem states convergence in an $L^2$-sense the convergence in the weak sense (\ref{metrdef}) is easy to conclude.\\[2mm]
\textbf{Corollary.}\it\ \ With the conditions of the theorem we have\rm
\bea\label{stat2}
\limn E\left[\sur\left(\sup_{t\in [0,T^* ]} D(\snr (t),S_r (.,t))\right)\right]=0.
\eea

\section{Proof of the convergence result}
  \subsection{Preliminaries}
The following lemma is useful in central estimates of section IV.\\[1mm]
\textbf{Lemma.}\it\ \ With $f\in L^2 \cap C_b^1 (\rel^d ,\rel)$, $\nab f\in L^2 (\rel^d ,\rel^d )$ we have\rm
\bea\label{lemma1}
\| f-f*W_N \|_2^2 \leq C\alpha_N^{-2} \|\nab f\|_2^2 \quad\mbox{and}\quad\| f-f*W_N \|_\infty \leq C\alpha_N^{-1} \|\nab f\|_\infty.
\eea\it
An analogous estimate is true, if $W_N$ and $\alpha_N$ are replaced by $\hw_N$ and $\han$.\\
For any finite, positive measure $\mu$ on $\rel^d$ and with $U_N (\xbo):=|\xbo |\hw_N(\xbo )$ there is\rm
\bea\label{lemma2}
\|\mu *U_N \|_2^2 \leq C\han^{2\epsilon -2} \|\mu *\hw_N \|_2^2 +\langle\mu ,1\rangle^2 \exp (-C'\han^\epsilon )\qquad\mbox{\it for all \rm}\epsilon>0.
\eea\it
For any finite, signed measure $\mu$ on $\rel^d$ it is\rm
\bea\label{lemma3}
\|\mu *W_N \|_2^2 \leq\|\mu *\hw_N\|_2^2 .
\eea
\textbf{Proof.} see \cite{oel1} sect.\ 4A,B and \cite{oel4} sect.\ 5B, or \cite{stefan}, sect.\ 4.3\\
In the proof of the lemma there is essentially made use of the conditions (C1) and (C2) on the interaction function $W_N$ and the kernel $\hw_N$. Due to the conservation of mass in the microscopic system (\ref{micrcons}) and with (\ref{ndef}) we get the following bound on the empirical processes,
\bea\label{sbound}
\sur\langle\snr(t),1\rangle =\oon\sum_{k\in M(N,t)} 1\leq {N\over N}=1.
\eea
We also use the following property without explicitly noting it for all suitable $f$ and $g$, such that the expressions are well defined:
\bea\label{symm}
\langle f,g*W_N \rangle =\langle f*W_N ,g\rangle\quad\mbox{and}\quad\langle f,g*\hw_N \rangle =\langle f*\hw_N ,g\rangle ,
\eea
because $W_1$ is symmetric according to (C2). Throughout this chapter $C$, $C'$ etc.\ denote suitably chosen constants, whose value can vary from line to line.

  \subsection{Proof of the theorem}
To prove statement (\ref{stat1}) we first look at the time evolution of the quantity
\bea\label{start}
\|\dnr (.,t)\|_2^2 =\| h_{N,r} (.,t)\|_2^2 -2\langle h_{N,r} (.,t),s_r (.,t)\rangle +\| s_r (.,t)\|_2^2 .
\eea
The dynamics of the first two terms is obtained analogous to (\ref{micr4}) using (\ref{micr1}), (\ref{jump}), (\ref{symm}) and It\^o's formula:
\bea
\|\hnr (.,t)\|_2^2 &=&{1\over N^2} \sum_{k,l\in M(N,r,t)}\big(\hw_N *\hw_N \big) (\xnk(t)-\xnl(t))\nonumber\\
\langle\hnr (.,t),s_r (.,t)\rangle &=&\oon\sukrt \big( s_r (.,t)*\hw_N \big) (\xnk(t))\nonumber
\eea
We just have to replace the test function $f$ in (\ref{micr4}) by $\snr (t)*(\hw_N *\hw_N )$ resp.\ $s_r (.,t)*\hw_N$. The expansion of the third term in (\ref{start}) follows from the macroscopic equation (MA):
\bea
\lefteqn{\|s_r (.,t)\|_2^2 =\|s_r (.,0)\|_2^2 +\int_0^t ds\Big\langle s_r (.,s),-2\,\nab\cdot \big(\vbo_r (.,s)s_r (.,s)\big) +\sigma_r^2 \Delta s_r (.,s)\Big\rangle}\nonumber\\
& &\quad +\sum_{q=1}^R \int_0^t ds\left\langle s_r (.,s),-2s_r (.,s)\,\harq(.,s)\, s_q (.,s)+s_q (.,s)\sul\heqlr\,\haql(.,s)\, s_l (.,s)\right\rangle\nonumber
\eea
Combining the parts suitably by using (\ref{micrcons}) to express $\heqlr$ in terms of $\eqlr$ we get:
\bea\label{dev}
\lefteqn{\|\dnr (.,t)\|_2^2 =\|\dnr (.,0)\|_2^2 }\nonumber\\
& &+\sigma_r^2 \int_0^t ds\ \big\langle\dnr (.,s),\Delta\dnr (.,s)\big\rangle\nonumber\\
& &+2\int_0^t ds\Big(\big\langle S_{N,r} (s),\nab (\dnr (.,s)*\hw_N )\cdot\vbo_r (.,s)\big\rangle -\big\langle s_r (.,s),\nab\dnr (.,s)\cdot\vbo_r (.,s)\big\rangle\Big)\nonumber\\
& &-2\sum_{q=1}^R \int_0^t ds\ \Big(\big\langle S_{N,r} (s),(\dnr (.,s)*\hw_N )\, \ta_{N,rq} (.,s)\big\rangle\nonumber\\
& &\qquad\qquad\qquad\quad -\big\langle s_r (.,s),\dnr (.,s)\,\harq(.,s)\, s_q (.,s)\big\rangle\Big)\nonumber\\
& &+2\sum_{q,l=1}^R \int_0^t ds\ \eqlr\Big(\big\langle S_{N,q} (s),(\dnr (.,s)*\hw_N )\, a_{N,ql}(.,s)\big\rangle\nonumber\\
& &\qquad\qquad\qquad\qquad\quad -\big\langle s_q (.,s),\dnr (.,s)\,\haql(.,s)\, s_l (.,s)\big\rangle\Big)\nonumber\\
& &+2\int_0^t {\sigma_r \over N}\sukrs\nab (\dnr (.,s)*\hw_N )(\xnk(s))\cdot d\bk (s)\nonumber\\
& &-2\sum_{q=1}^R \int_0^t \oon\sukrs (\dnr (.,s)*\hw_N )(\xnk(s))\left(\ta_{N,rq}^{*,k} (ds)-\ta_{N,rq} (\xnk(s),s)\, ds\right)\nonumber\\
& &+2\!\!\sum_{q,l=1}^R \!\int_0^t \!\oon\!\sukqs\!\!\!\!\!\! (\dnr (.,s)*\hw_N )(\xnk(s))\,\eqlr\!\left(\ta_{N,ql}^{*,k} (ds)-\ta_{N,ql} (\xnk(s),s)\, ds\right)\nonumber\\
& &-\sigma_r^2 {\Delta (\hw_N *\hw_N )(\vec{0})\over N} \int_0^t ds\ \langle\snr (s),1\rangle\nonumber\\
& &+{(\hw_N *\hw_N )(\vec{0})\over N} \suq\int_0^t ds\ \big\langle\snr(s),\anrq (.,s)\big\rangle\nonumber\\
& &+{(\hw_N *\hw_N )(\vec{0})\over N} \sum_{q,l=1}^R \int_0^t ds\ \eqlr\big\langle\snq(s),\anql (.,s)\big\rangle\nonumber\\
& &+{(\hw_N *\hw_N )(\vec{0})\over N^2} \sum_{q=1}^R \int_0^t \sukrs\left(\ta_{N,rq}^{*,k} (ds)-\ta_{N,rq} (\xnk(s),s)\, ds\right)\nonumber\\
& &+{(\hw_N *\hw_N )(\vec{0})\over N^2} \sum_{q,l=1}^R \int_0^t \sukqs\eqlr\left(\ta_{N,ql}^{*,k} (ds)-\ta_{N,ql} (\xnk(s),s)\, ds\right) =\nonumber\\
&=&\|\dnr (.,0)\|_2^2+\int_0^t \Big( T_{N,r}^1 (s)+T_{N,r}^2 (s)+T_{N,r}^3 (s)+T_{N,r}^4 (s)\Big)\, ds\nonumber\\
& &\qquad\qquad\quad\ \ \, +2\left( \hm_{N,r}^1 (t) +\hm_{N,r}^{2a} (t)+\hm_{N,r}^{2b} (t)\right)\nonumber\\
& &\qquad\qquad\quad\ \ \, +T_{N,r}^0 (t)+T_{N,r}^{0a} (t)+T_{N,r}^{0b} (t)+\hm_{N,r}^{0a} (t)+\hm_{N,r}^{0b} (t)
\eea
The terms in the above sum are labeled line by line. $T_N^1$ derives from the diffusion due to Brownian motion, $T_N^2$ from the particle transport, $T_N^3$ from the loss and $T_N^4$ from the gain of particles due to interacions. The stochastic integrals $\hm_{N,r}^1$, $\hm_{N,r}^{2a}$ and $\hm_{N,r}^{2b}$ represent the fluctuations due to stochasticity in the free particle dynamics and the interaction. The remaining terms are corrections resulting from the expansion of $\|\hnr (.,t)\|_2^2$. With It\^o's formula and (\ref{micr1}) we get for $T_{N,r}^0$:
\bea
\lefteqn{{1\over N^2}\sum_{k,l\in M(N,r,t)} d(\hw_N *\hw_N )(\xnk(t)-\xnl(t))=}\nonumber\\
&=&{1\over N^2}\!\!\!\sum_{k,l\in M(N,r,t)\atop k\neq l}\Big(\!\nab (\hw_N *\hw_N )(\xnk(t)\! -\!\xnl(t))\, (d\xnk(t)-d\xnl(t))\nonumber\\
& &\qquad\qquad\qquad\quad +{\sigma_r^2 \over 2}\Delta (\hw_N *\hw_N )(\xnk(t)\! -\!\xnl(t))\, (dt+dt)\!\Big)\! =\nonumber\\
&=&2\big\langle\snr(t),\nab (\hnr(.,t)*\hw_N )\cdot \vbo_r (.,t)\big\rangle\, dt +{2\sigma_r\over N}\!\!\!\sukrt\!\!\!\!\!\nab\big(\hnr (.,t)*\hw_N \big)\cdot d\bk (t)\nonumber\\
& &+\sigma_r^2 \big\langle\hnr(.,t),\Delta\hnr(.,t)\big\rangle\, dt -{\sigma_r^2 \over N^2}\sum_{k\in M(N,r,t)}\Delta (\hw_N *\hw_N )(\vec{0})\, dt\nonumber
\eea
$T_{N,r}^{0a}$ and $T_{N,r}^{0b}$ can be derived analogously by considering the change of the sets $M(N,q,t)$, $q=1,\ldots ,R$. The fluctuations of these corrections are separated in stochastic integrals $\hm_{N,r}^{0a}$ and $\hm_{N,r}^{0b}$.

    \subsubsection*{Estimate of the correction terms:}
With (C1), (C2) and (\ref{sbound}) we get:
\bea
|T_{N,r}^0 (t)|&\leq& {t \sigma_r^2\over N}\left( \int_{\rel^d} |\hw_N (\xbo )\ \Delta \hw_N (\vec{0}-\xbo)|\, d^d x\right)\sup_{s\leq t}\langle\snr(s),1\rangle\leq\nonumber\\
&\leq& {C\over N} t\int_{\rel^d} N^{\hbeta (1+2/d)} \|\Delta W_1 \|_\infty |\hw_N (\xbo )|\, d^d x\leq Ct N^{\hbeta (1+2/d)-1} \nonumber
\eea
This term vanishes in the limit $N\rightarrow\infty$, because with (C1) it is $\hbeta <d/(d+2)$. Using also (\ref{adef}) and (\ref{micrcons}) we get completely analogous:
\bea
|T_{N,r}^{0a} (t)|+|T_{N,r}^{0b} (t)|\leq Ct N^{\hbeta -1}\nonumber
\eea
The stochastic integrals are estimated in section IV.D.

    \subsubsection*{Estimate of $T_{N,r}^1 (s)=\sigma_r^2 \langle\dnr (.,s),\Delta\dnr (.,s)\rangle$}
After partial integration we get with condition (C7): $T_{N,r}^1 (s)=-\sigma_r^2 \|\nab\dnr (.,s)\|_2^2 <0$. This term is negative and can be used to cancel positive contributions of the same kind arising in the estimates of $T_{N,r}^2$ and $\hm_{N,r}^1$.

    \subsubsection*{Estimate of $T_{N,r}^2 (s)\!\! =\!\! 2\Big(\!\big\langle S_{N,r} (s),\!\nab (\dnr (.,s)\! *\!\hw_N )\!\cdot\!\vbo_r (.,s)\big\rangle\! -\!\big\langle s_r (.,s),\!\nab\dnr (.,s)\!\cdot\!\vbo_r (.,s)\big\rangle\!\Big)$}
To contract the two brackets we make the following replacement:
\bea
\big\langle S_{N,r} (s),(\nab \dnr (.,s)*\hw_N )\cdot\vbo_r (.,s)\big\rangle =\big\langle h_{N,r} (.,s),\nab \dnr (.,s)\cdot\vbo_r (.,s)\big\rangle +R_{N,r}^2 (s),\nonumber
\eea
where the correction term is estimated using (C6):
\bea
\lefteqn{|R_{N,r}^2 (s)|=\left|\left\langle S_{N,r} (s),\int_{\rel^d} d^d u\, \hw_N (\ubo) \nab \dnr (.-\ubo ,s)\cdot \big(\vbo_r (.,s)-\vbo_r (.-\ubo ,s)\big)\right\rangle\right| \leq}\nonumber\\
& &\qquad\leq\left\langle S_{N,r} (s),\int_{\rel^d} d^d u\, \hw_N (\ubo) |\ubo| |\nab \dnr (.-\ubo,s)| \right\rangle \|\nab\vbo_r (.,s)\|_\infty \leq\nonumber\\
& &\qquad\leq C \big\langle S_{N,r} (s)*U_N ,|\nab \dnr (.,s)| \big\rangle\leq C \left(\tilde{C} \| S_{N,r} (s)*U_N \|_2^2 +{1\over\tilde{C}} \| \nab \dnr (.,s)\|_2^2 \right) ,\nonumber
\eea
with $U_N (\xbo ):=|\xbo |\hw_N (\xbo )$. This is true for all $\tilde{C}>0$ using Cauchy's inequality. With the second statement of the lemma (\ref{lemma2}) and (\ref{sbound}) we have for all $\epsilon >0$
\bea
|R_{N,r}^2 (s)|&\leq&C\tilde{C} \Big(\han^{2\epsilon -2} \| h_{N,r} (.,s)\|_2^2 + \langle S_{N,r} (s),1\rangle^2 e^{-C' \han^\epsilon }\Big) +{C'' \over \tilde{C}} \| \nab \dnr (.,s)\|_2^2 \leq\nonumber\\
&\leq&C\tilde{C} \Big( \han^{2\epsilon -2} (\| \dnr (.,s)\|_2^2 + C') + e^{-C'' \han^\epsilon }\Big) +{C''' \over \tilde{C}} \| \nab \dnr (.,s)\|_2^2 ,\nonumber
\eea
using (C3), because with the triangle inequality it is
\bea
\|\hnr(.,s)\|_2^2 \leq \big(\|\dnr(.,s)\|_2 +\| s_r (.,s)\|_2 \big)^2 \leq 2\|\dnr(.,s)\|_2^2 +2\| s_r (.,s)\|_2^2 .\nonumber
\eea
$\tilde{C}$ is chosen after the estimate of $\hm_{N,r}^1$, so that the term arising there and ${C'''\over \tilde{C}} \| \nab \dnr (.,s)\|_2^2$ cancels with the negative contribution from the estimate of $T_{N,r}^1$. Choosing $\epsilon =\halb$ the constant terms in the above estimate of $|R_{N,r}^2 (s)|$ and the prefactor of $\| \dnr (.,s)\|_2^2$ vanish in the limit $N\rightarrow\infty$. Now we can write
\bea
T_{N,r}^2 (s)=2\langle\dnr (.,s),(\nab\dnr (.,s))\cdot \vbo_r (.,s)\rangle +2 R_{N,r}^2 (s).\nonumber
\eea
With (C6) and the estimate
\bea
\left|\big\langle \dnr (.,s),(\nab \dnr (.,s))\cdot\vbo_r (.,s)\big\rangle\right| &=&\halb \left|\big\langle \nab d_{N,r}^2 (.,s),\vbo_r (.,s)\big\rangle\right| =\nonumber\\
=\left|-\halb\big\langle d_{N,r}^2 (.,s),\nab\cdot\vbo_r (.,s)\big\rangle\right| &\leq&\halb\|\nab\cdot\vbo_r (.,s)\|_\infty \|\dnr (.,s)\|_2^2\nonumber
\eea
we get after a suitable arrangement of constants and terms:
\bea
|T_{N,r}^2 (s)|\leq C\big( 1 +\tilde{C}\big)\|\dnr (.,s)\|_2^2 +{C' \over\tilde{C}} \|\nab \dnr (.,s)\|_2^2 +\tilde{C}\ O\big( N^{-\hbeta/d}\big) ,\nonumber
\eea
because the $N$-dependent prefactors of $\|\dnr (.,s)\|_2^2$ vanish monotonically with $N\rightarrow\infty$ and $N^{-\hbeta/d}$ is the leading order in $N$ of all constant terms.

    \subsubsection*{Estimate of $T_{N,r}^3 (s)=-2\sum_{q=1}^R \Big(\big\langle S_{N,r} (s),(\dnr (.,s)*\hw_N )\, \ta_{N,rq} (.,s)\big\rangle$\\
\hspace*{55mm}$-\big\langle s_r (.,s),\dnr (.,s)\,\harq(.,s)\, s_q (.,s)\big\rangle\Big)$}
First we make the same substitution as before:
\bea
\big\langle S_{N,r} (s),(\dnr (.,s)*\hw_N )\,\ta_{N,rq} (.,s)\big\rangle =\big\langle h_{N,r} (.,s), \dnr (.,s)\,\ta_{N,rq} (.,s)\big\rangle +R_{N,rq}^{3a} (s),\nonumber
\eea
and we can estimate the correction term analogously to $R_{N,r}^2$ using (\ref{adef}), (C1), (C2), (C4), (\ref{lemma2}) and (\ref{sbound}):
\bea
|R_{N,rq}^{3a} (s)|&=&\left|\left\langle S_{N,r} (s),\!\int_{\rel^d}\!\! d^d u\, \hw_N (\ubo)\ \dnr (.-\ubo,s) \big(\ta_{N,rq} (.,s)-\ta_{N,rq} (.-\ubo,s)\big)\right\rangle\right|\leq\nonumber\\
&\leq&\left\langle S_{N,r} (s),\int_{\rel^d} d^d u\,\hw_N (\ubo) |\ubo| |\dnr (.-\ubo,s)|\right\rangle\nonumber\\
& &\Big(\|\harq(.,s)\|_\infty \|\nab (\snq(s)\! *\! W_N )\|_\infty\! +\!\|\nab\harq(.,s)\|_\infty\|\snq(s)\! *\! W_N \|_\infty \Big)\leq\nonumber\\
&\leq&\big\langle S_{N,r} (s)*U_N ,|\dnr (.,s)|\big\rangle\ C\ \big(\alpha_N^{d+1} +\alpha_N^d \big)\big(\| W_1 \|_\infty +\|\nab W_1 \|_\infty\big)
\leq\nonumber\\
&\leq&C \alpha_N^{2d+2} \Big(\han^{2\epsilon -2} (\| \dnr (.,s)\|_2^2 + C' )+ e^{-C'' \han^\epsilon }\Big) + \| \dnr (.,s)\|_2^2 \nonumber
\eea
For the first term in the above estimate to vanish in the limit $N\rightarrow\infty$, we choose $\epsilon >0$ so that $\gamma :=(\beta (2d+2)+\hbeta (2\epsilon -2))/d<0$. That means $0<\epsilon <{1\over\hbeta}(\hbeta -\beta (d+1))$, which is possible due to condition (C1) on the scaling parameters $\beta$ and $\hbeta$. Now we look at the remaining term in $T_{N,r}^3$:
\bea
-2\sum_{q=1}^R \Big(\big\langle h_{N,r} (.,s),\dnr (.,s)\,\ta_{N,rq} (.,s)\big\rangle -\big\langle s_r (.,s),\dnr (.,s)\,\harq(.,s)\ s_q (.,s)\big\rangle\Big) \nonumber
\eea
To contract the two brackets we have to compare the microscopic and macroscopic interaction rates:
\bea\label{inter}
\big\langle h_{N,r} (.,s),\dnr (.,s)\,\ta_{N,rq} (.,s)\big\rangle -\big\langle s_r (.,s),\dnr (.,s)\,\harq (.,s)\, s_q (.,s)\big\rangle = \nonumber\\
=\big\langle \dnr (.,s),\dnr (.,s)\,\ta_{N,rq} (.,s)\big\rangle +R_{N,rq}^{3b} (s),\nonumber
\eea
with the correction term
\bea
R_{N,rq}^{3b} (s)=\left\langle s_r (.,s),\dnr (.,s) \big(\ta_{N,rq} (.,s) -\harq (.,s)\, s_q (.,s)\big)\right\rangle .\nonumber
\eea
With the definition of $\ta_{N,rq} (.,s)$ in (\ref{adef}) and (C5) we get the following estimate:
\bea
|R_{N,rq}^{3b} (s)|&\leq&\!\left\langle\! s_r (.,s),\!|\dnr (.,s)| |\harq (.,s)|\!\left(\! |s_q (.,s)-\snq(s)\! *\! W_N |\! +\delta_{r,q} \oon W_N (\vec{0})\!\right)\!\right\rangle\!\leq\nonumber\\
&\leq&\|s_r (.,s)\,\harq(.,s)\|_\infty \|\dnr (.,s)\|_2 \| s_q (.,s)-\snq(s)*W_N \|_2 \nonumber\\
& &+\delta_{r,q} N^{\beta -1} W_1 (\vec{0}) \left(\|\harq(.,s)\|_\infty^2 \| s_r (.,s)\|_2^2 +\|\dnr(.,s)\|_2^2 \right)\nonumber
\eea
It is $\|s_q (.,s)-\snq(s)*W_N \|_2 \leq \|\big( s_q (.,s)-\snq (s)\big) *W_N \|_2+\|s_q (.,s)-s_q (.,s)*W_N \|_2$\linebreak and with the third statement (\ref{lemma3}) of the lemma
\bea
\left\|\big( s_q (.,s)-\snq (s)\big) *W_N \right\|_2 &\leq&\left\|\big( s_q (.,s)-\snq (s)\big) *\hw_N \right\|_2 \leq\nonumber\\
&\leq&\|\dnq (.,s)\|_2 +\|s_q (.,s)-s_q (.,s)*\hw_N \|_2 .\nonumber
\eea
Therefore we get with the first statement (\ref{lemma1}) and (C3)
\bea
|R_{N,rq}^{3b} (s)|&\leq& C \Big( \|\dnr (.,s)\|_2^2 +\|\dnq (.,s)\|_2^2 +C'(\alpha_N^{-2} +\han^{-2} )\|\nab s_q (.,s)\|_2^2 \Big)\nonumber\\
& &+\delta_{r,q} C'' N^{\beta -1} \left( 1+\|\dnr(.,s)\|_2^2 \right) .\nonumber
\eea
After treating the correction terms we get for the main contribution using (\ref{adef}):
\bea
|\langle \dnr (.,s),\dnr (.,s)\, \ta_{N,rq} (.,s)\rangle |\leq C_a \|\dnr (.,s)\|_2^2\nonumber
\eea
Arranging all terms analogously to $T_{N,r}^2$ we finally have the estimate
\bea
|T_{N,r}^3 (s)|\leq C\suq\|\dnq (.,s)\|_2^2 +O\left( N^\gamma\right)\ ,\qquad\mbox{with }\gamma <0.\nonumber
\eea

    \subsubsection*{Estimate of $T_{N,r}^4 (s)=2\sum_{q,l=1}^R \eqlr\Big(\big\langle S_{N,q} (s),(\dnr (.,s)*\hw_N )\, a_{N,ql}(.,s)\big\rangle$\\
\hspace*{60mm}$-\big\langle s_q (.,s),\dnr (.,s)\,\haql(.,s)\, s_l (.,s)\big\rangle\Big)$}
Obviously, the estimate of this term is completely analogous to the one of $T_N^3$ with the same result except for different constants,
\bea
|T_{N,r}^4 (s)|\leq C\suq\|\dnq (.,s)\|_2^2 +O\left( N^\gamma\right)\ ,\qquad\mbox{with }\gamma <0.\nonumber
\eea\\
Putting all the estimates together, taking the absolute value, the supremum over all $t\in [0,T]$ for some $T<T^*$ (see (C3)) and the sum $\sur$ on both sides of equation (\ref{dev}) we arrive at:
\bea\label{devesup}
\lefteqn{\sur\left(\sup_{t\in [0,T]} \| \dnr (.,t)\|_2^2 + \left(\sigma_r^2 -{C \over\tilde{C}}\right)\int_0^T ds\,\| \nab \dnr (.,s)\|_2^2 \right)\leq }\nonumber\\
&\leq&\!\!\sur \left(\| \dnr (.,0)\|_2^2 +T\ C'\sup_{t\in [0,T]} \| \dnr (.,t)\|_2^2 \right)  +(1+\tilde{C} +T)\, O\left( N^\gamma \right)\nonumber\\
& &\!\!\!\! +\sur\sup_{t\in [0,T]} \left( 2|\hm_{N,r}^1 (t)|+2|\hm_{N,r}^{2a} (t)|+2|\hm_{N,r}^{2b} (t)|+|\hm_{N,r}^{0a} (t)|+|\hm_{N,r}^{0b} (t)|\right)
\eea
The choice of $\tilde{C}$ below ensures the positivity of all occuring terms. The leading order of all constant terms, including the estimated correction terms, is characterized by $\gamma <0$ defined in the estimate of $T_{N,r}^3$. Taking the expectation on both sides of (\ref{devesup}) we can use the estimates (\ref{mres}) and (\ref{mres2}) of the stochastic integrals in section IV.D and get:
\bea\label{mdevesup}
E\left[\sur\left(\sup_{t\in [0,T]} \|\dnr (.,t)\|_2^2 +\left(\sigma_r^2 -C N^{(\hbeta -1)/2} -{C' \over\tilde{C}}\right)\int_0^T ds\| \nab \dnr (.,s)\|_2^2 \right)\right] \leq\nonumber\\
\leq E\left[\sur\left(\| \dnr (.,0)\|_2^2 +T\, C''\sup_{t\in [0,T]} (\| \dnr (.,t)\|_2^2 )\right)\right] +(1+\tilde{C} +T)\ O\left( N^\gamma \right) ,\quad
\eea
where the leading order of constant terms remains unchanged with (C1). Now we can choose $\tilde{C}$, occuring in the estimate of $T_{N,r}^2$, and $N_0 \in\nat$ large enough, so that the prefactor of\linebreak $\int_0^T ds\,\| \nab \dnr (.,s)\|_2^2$ on the lefthand side of (\ref{mdevesup}) is positive for all $N>N_0$ (see (C7)). Consequently this term can be neglected and after a short rearrangement we have for all $N>N_0$ and $0<T<1/C''$, where $C''$ is taken from (\ref{mdevesup}):
\bea\label{fertig}
E\left[\sur\sup_{t\in [0,T]} \|\dnr (.,t)\|_2^2 \right]\!\!\leq\! {1\over 1-T\, C''}\, E\left[\sur\| \dnr (.,0)\|_2^2 \right] +(1+T)\ O\left( N^\gamma \right)
\eea
Taking the limit $N\rightarrow\infty$ on both sides, the constant terms vanish as $\gamma<0$ and convergence at $t=0$ is given in the theorem. This finally proves statement (\ref{stat1}) for all $t\in [0,T]$, but the above constraint on the time interval is not essential. At $t=T$ all the conditions for the theorem are fulfilled as long as $T<T^*$, which enables us to apply the proof again with the same constants arising. The length of the next time interval is subject to the same constraint and after a finite number of repetitions this proves the theorem.

  \subsection{Proof of the corollary}
We have the following estimate:
\bea
|\langle\snr (t)-s_r (.,t),f\rangle|&\leq& |\langle\snr (t)-\hnr(.,t),f\rangle |+|\langle\dnr (.,t),f\rangle | \leq\nonumber\\
&\leq& |\langle\snr (t),f-f*\hw_N \rangle |+ \|f\|_2 \|\dnr (.,t)\|_2 \leq\nonumber\\
&\leq& \langle\snr (t),1\rangle\| f-f*\hw_N \|_\infty +\|f\|_2 \|\dnr (.,t)\|_2 \nonumber
\eea
Therefore we get for all $f\in\Omega_D$ (see \ref{metrdef}), using statement (\ref{lemma1}) of the lemma and (\ref{sbound}): $|\langle\snr (t)-s_r (.,t),f\rangle|\leq C(\|\nab f\|_\infty \han^{-1} +\|\dnr (.,t)\|_2 )$. Hence we have for the metric\nolinebreak $\ D$
\bea\label{ungl}
|D(\snr (t),S_r (t))|\leq C(\han^{-1} +\|\dnr (.,t)\|_2 )\quad\mbox{for all }r=1,\ldots ,R.
\eea
It is straightforward to see that the statement (\ref{stat1}) of the theorem implies
\bea
\limn E\left[\sur\sup_{t\in [0,T^* ]} \|\dnr (.,t)\|_2 \right] =0.\nonumber
\eea
Using this, equation (\ref{ungl}) immediately gives
\bea
\limn E\left[\sur\left(\sup_{t\in [0,T^* ]} D(\snr (t),S_r (t))\right)\right]\leq C\limn E\left[\sur\sup_{t\in [0,T^* ]} \|\dnr (.,t)\|_2 \right]=0,\nonumber
\eea
which proves the statement (\ref{stat2}) of the corollary.

  \subsection{Estimate of the stochastic integrals}
In the following we estimate the stochastic integrals occuring in (MI) and (\ref{dev}). First we show for every fixed system size $N\in\nat$ that the integrands and integrators fulfill the necessary regularity conditions so that standard techniques of stochastic integration taken from \cite{ikeda,karatzas} can be applied. This part is kept short and can be read in more detail in \cite{stefan} (sect.\ 3.1.1, appendices B and C), following the work in \cite{oel1}. According to those results we find an estimate uniformly in $N$ using Doob's inequality \cite{karatzas}. First we consider the terms
\bea\label{integrals}
M_{N,r}^1 (f,t)&=&{\sigma_r \over N}\sum_{k=1}^\infty \sum_{i=1}^d \int_0^t I_{N,r,k,i}^1 (f,s)\ dB_i^k (s)\nonumber\\
& &\mbox{where }I_{N,r,k,i}^1 (f,t):=\chi_{M(N,r,t)} (k)\,\partial_{x_i} f(\xnk(t))\nonumber\\
M_{N,r}^{2a} (f,t)&=&-{1\over N}\suq\sum_{k=1}^\infty \int_0^t I_{N,r,k}^{2a} (f,s)\ dP_{N,rq}^k (s)\nonumber\\
& &\mbox{where }I_{N,r,k}^{2a} (f,t):=\chi_{M(N,r,t)} (k)\, f(\xnk(t))\nonumber\\
& &\mbox{and  }P_{N,rq}^k (t):=\anrqs (t)-\int_0^t \chi_{M(N,r,s)} (k)\,\anrq (\xnk(s),s)\ ds\qquad
\eea
occuring in (MI). $M_{N,r}^{2b}$ can be handled analogously to $M_{N,r}^{2a}$. The integrators $B_i^k$ and $P_{N,rq}^k$ in (\ref{integrals}) are Brownian motions and time-inhomogeneous compensated Poisson processes (\ref{jump}) with bounded rates (\ref{adef}). Therefore they are square integrable martingales with respect to the filtration $(\mathcal{F}_t )_{t\geq 0}$, mentioned after equation (\ref{micr1}). The integrands $I_{N,r,k,i}^1$ and $I_{N,r,k}^{2a}$ can be replaced by their left continuous versions as we have
\bea
\lefteqn{E\left[\int_0^t |I_{N,r,k,i}^1 (f,s)-I_{N,r,k,i}^1 (f,s-)|^2 d\langle B_i^k \rangle (s)\right] =}\nonumber\\
& &\qquad =E\left[\int_0^t |I_{N,r,k}^{2a} (f,s)-I_{N,r,k}^{2a} (f,s-)|^2 d\langle P^k_{N,rq} \rangle (s)\right] =0\ ,\nonumber
\eea
because the quadratic variational processes $\langle B_i^k \rangle (s)=s$ and $\langle P^k_{N,rq} \rangle (s)=\int_0^s \chi_{M(N,r,u)}(k)$\linebreak $a_{N,rq} (\xnk(u),u)\, du$ are continuous in $s$. On Poisson processes with time dependent rates see e.g.\ \cite{ikeda}, theorem II.3.1 on page 60. Therefore we can assume the integrands to be predictable processes with respect to $(\mathcal{F}_t )_{t\geq 0}$. For every fixed system size $N\in\nat$ it is
\bea
E\left[\int_0^t |I_{N,r,k,i}^1 (f,s)|^2 d\langle B_i^k \rangle (s)\right]\ ,\ E\left[\int_0^t |I_{N,r,k}^{2a} (f,s)|^2 d\langle P^k_{N,rq} \rangle (s)\right] <\infty ,
\eea
using $f\in C_b^2 (\rel^d ,\rel)$ and estimates analogous to (\ref{m1}) and (\ref{m2}) below. Hence $M_{N,r}^1 (f,t)$ and\linebreak $M_{N,r}^{2a} (f,t)$ are martingales with respect to $(\mathcal{F}_t )_{t\geq 0}$. On the other hand we have for all $r=1,\ldots ,R$, $t\geq 0$ and $N\in\nat$:
\bea\label{m1}
\lefteqn{E\left[ |M_{N,r}^1 (f,t)|^2 \right] =}\nonumber\\
& &\quad =E\left[{\sigma_r^2 \over N^2}\!\left(\!\int_0^t \!\!\sukrs\!\!\!\!\!\!\!\nab f(\xnk(s))\!\cdot\! d\bk (s)\!\right)\!\!\left(\!\int_0^t \!\!\sum_{l\in M(N,r,s)}\!\!\!\!\!\!\!\nab f(\xnl(s))\!\cdot\! d\vec{B}^l (s)\!\right)\!\right]\nonumber\\
& &\quad ={\sigma_r^2 \over N^2} E\left[\int_0^t \sum_{k,l\in M(N,r,s)} \sum_{i,j=1}^d \partial_{x_i} f(\xnk(s))\ \partial_{x_j} f(\xnl(s))\ d\langle B_i^k ,B_j^l \rangle (s)\right] =\nonumber\\
& &\quad ={\sigma_r^2 \over N^2} E\left[\int_0^t \sukrs|\nab f(\xnk(s))|^2 ds\right]={\sigma_r^2 \over N} \int_0^t E[\langle\snr (s),|\nab f|^2 \rangle ]\ ds\leq\nonumber\\
& &\quad\leq {\sigma_r^2 \over N} \|\nab f\|_\infty^2 \int_0^t E[\langle\snr (s),1\rangle ]\ ds\leq C\, {\sigma_r^2 \over N} \|\nab f\|_\infty^2 t<\infty
\eea
\bea\label{m2}
\lefteqn{E\left[ |M_{N,r}^{2a} (f,t)|^2 \right] =}\nonumber\\
&=&\!\! {1\over N^2} \sum_{p,q=1}^R E\left[\int_0^t \sum_{k,l\in M(N,r,s)} f(\xnk(s))f(\xnl(s))\, d\langle P_{N,rp}^k ,P_{N,rq}^l \rangle\right] =\nonumber\\
&=&\!\! {1\over N^2} \suq E\left[\int_0^t \sukrs f^2 (\xnk(s))\ a_{N,rq} (\xnk(s),s)\, ds\right]\leq\nonumber\\
&\leq&\!\!\oon\|f\|_\infty^2 \!\left(\!\sup_{s\in [0,t]} \suq \| a_{N,rq} (.,s)\|_\infty\! \right)\!\int_0^t \!\!\!\! E\left[\langle\snr(s),1\rangle\right] ds\leq{RC_a \over N} \|f\|_\infty^2 t<\infty
\eea
This is true with (\ref{sbound}), $f\in C_b^2 (\rel^d ,\rel)$, (C2) and (C4). After an estimate of $M_{N,r}^{2b}$ analogous to $M_{N,r}^{2a}$ we can apply Doob's inequality and get for $i=1,2a,2b$, $r=1,\ldots ,R$ and all $T>0$:
\bea\label{doobm}
E\left[\sup_{t\in [0,T]} |M_{N,r}^i (f,t)|\right]^2 \leq E\left[\sup_{t\in [0,T]} |M_{N,r}^i (f,t)|^2\right] \leq 4E\big[ |M_{N,r}^i (f,T)|^2 \big] ,
\eea
so the result (\ref{erg1}) follows from (\ref{m1}) to (\ref{doobm}).\\
For the terms occuring in (\ref{dev}) we have $\hm_{N,r}^i (t)=M_{N,r}^i (\dnr(.,t)*\hw_N ,t)$ for $i=1,2a,2b$ and therefore a similar reasoning applies in this case. To get an estimate in the limit $N\rightarrow\infty$ we note that with (C1) and (C2) $\hw_N (\xbo -.)$ is a probability density on $\rel^d$ for all $\xbo\in\rel^d$. Hence we have for all $f:\rel^d \rightarrow\rel$, where $E_{\hw_N}$ denotes the corresponding expectation value:
\bea
|(f*\hw_N )(\xbo)|^2 =\left|\int_{\rel^d} \hw_N (\xbo-\ybo)f(\ybo)\, d^d y\right|^2 =E_{\hw_N} [f]^2 \leq E_{\hw_N} [f^2]=(|f|^2 *\hw_N )(\xbo)\nonumber
\eea
We also have with (C1), (C2) and (\ref{sbound}) for all $t\in [0,T^* ]$ and $r=1,\ldots ,R$:
\bea
\|\hnr(.,t)\|_\infty \leq\han^d \oon\!\!\sum_{k\in M(N,r,t)}\!\! \| W_1 (\han (.-\xnk(t)))\|_\infty \leq C N^\hbeta \langle\snr (t),1\rangle\leq C N^\hbeta \nonumber
\eea
With the last two considerations we get for all $T\in [0,T^* ]$, starting with the third line of (\ref{m1}) and second line of (\ref{m2}) and using (\ref{symm}),
\bea
\lefteqn{E\left[ |\hm_{N,r}^1 (T)|^2 \right]\leq{\sigma_r^2 \over N} E\left[\int_0^t \langle\snr(s),|\nab (\dnr (.,s)*\hw_N )|^2 \rangle\, ds\right]\leq}\nonumber\\
& &\quad\leq {\sigma_r^2 \over N}E\left[\int_0^T \!\!\langle\hnr (.,s),|\nab\dnr (.,s)|^2 \rangle\, ds\right]\!\leq CN^{\hbeta -1} \sigma_r^2 E\left[\int_0^T \!\!\|\nab\dnr (.,s)\|_2^2 \, ds\right] ,\nonumber
\eea
\bea
\lefteqn{E\left[ |\hm_{N,r}^{2a} (T)|^2 \right]\leq{1\over N} \suq E\left[\int_0^t \langle\snr(s), |\dnr(.,s)*\hw_N |^2 \ a_{N,rq} (.,s)\, ds\right]\leq}\nonumber\\
& &\quad\leq{RC_a \over N} E\left[\int_0^T \langle\hnr(.,s),|\dnr(.,s)|^2 \rangle ds\right]\leq CN^{\hbeta -1} E\left[\int_0^T \|\dnr (.,s)\|_2^2 \, ds\right] .\nonumber
\eea
Using Doob's inequality in the form (\ref{doobm}) we get for all $T\in [0,T^* ]$:
\bea\label{mres}
E\left[\sur\sup_{t\in [0,T]} |\hm_{N,r}^1 (t)| \right]&\leq& CN^{(\hbeta -1)/2} \left( 1+\sur E\left[\int_0^T \|\nab\dnr (.,s)\|_2^2 \, ds\right]\right)\nonumber\\
E\left[\sur\sup_{t\in [0,T]} |\hm_{N,r}^{2a} (t)|\right]&\leq& C\left( N^{\hbeta -1} +T\sur E\left[\sup_{t\in [0,T]} \|\dnr(.,t)\|_2^2 \right]\right) .
\eea
With (\ref{micrcons}) the analogous estimate of $\hm_N^{2b}$ yields the same result as for $\hm_N^{2a}$ except for the constant $C$. The correction terms $\hm_{N,r}^{0a}$ and $\hm_{N,r}^{0b}$ are much easier to handle and with (C1), (C2), (\ref{m2}) and (\ref{micrcons}) we get
\bea\label{mres2}
E\left[\sur\left(\sup_{t\in [0,T]} |\hm_{N,r}^{0a} (t)|+\sup_{t\in [0,T]} |\hm_{N,r}^{0b} (t)|\right)\right] \leq CT N^{\hbeta -2}
\eea

\section{Discussion}
  \subsection{Connection to Astrophysics}
In the astrophysical context studied in \cite{jones,yorke1} the particles are dust grains in a star forming cloud of hydrogen gas. Depending on the grain size distribution these particles determine the opacity of the gas cloud and influence thermodynamic properties by emission and absorption of heat radiation. They also interact with the gas via hydrodynamic drag and influence chemical reactions via catalysis. A realistic expression for the collision rate of two ball shaped grains with radii $l_r$ and $l_q$ is $\harq (\xbo ,t)\sim (l_r +l_q )^{d-1} g\big( |\vbo_r (\xbo ,t)-\vbo_q (\xbo ,t)|\big)$, proportional to the cross section and depending on the relative velocity of the two particles (cf.\ \cite{yorke1}). The material coefficient $\herql$ should also depend on the relative velocity of the two collision partners. If the latter is high, shattering collisions with several outgoing particles of masses $m_l <m_r ,m_q$ are prefered. If it is low, the colliding particles are more likely to coagulate, so that there is one outgoing particle of mass $m_l =m_r +m_q$. For a precise form of this function from empirical data for different grain materials we refer to \cite{jones} (table 1) and references therein.\\
Of special interest in \cite{jones,yorke1} is the situation after a shock with very high gas temperatures and particle densities. In this regime shattering is the dominant process, justifying (C4), and it seems natural to assume that all particles are Brownian (C7). The stochasticity in the interaction coming from regularity conditions on $W_1$ in (M3) also appears to be reasonable. The effective field interaction introduced in (M4) is a simplification we have to make in order to include the space dependence in $\harq$, which cannot be included in a model with pair interactions (see sect.\ V.B). The cut-off of the interaction rate in (M4) prevents a divergence due to high particle concentrations and seems natural, as for real grains the density is limited due to the positive particle volume. As the focus of this paper is on the interaction and not on transport terms in (MA), we left out the complicated hydrodynamic drag interaction between particles and gas in (M2). We note that all realistic features are covered by our derivation, except for the velocity dependence of $\herql$, which can be included in a direct generalization explained in the next subsection.\\
The moderate interaction scaling is technically important, as seen in the proof, and can also be interpreted on a physical level. In the limit $N\rightarrow\infty$ the scaling of the mean distance between particles is given by $\Lambda_N \sim N^{-1/d}$ and due to (\ref{rate}) the interaction radius scales like $r_N \sim N^{-\beta /d}$ with $\beta <1$. It vanishes more slowly with $N$ than $\Lambda_N$ does and the number of interaction partners of a particle diverges in the limit $N\rightarrow\infty$. In contrast to the hydrodynamic scaling with $\Lambda_N \sim r_N$ this introduces a self-averaging effect and the influence of the interaction partners is determined by the local particle density. In the astrophysical gas cloud after a shock the cross section for a collision interaction can be thought of being effectively enlarged by strong Brownian motion of the particles. Together with high particle densities this leads to a large number of interaction partners and justifies the moderate scaling in our microscopic model (M3), whereas the scaling in the Boltzmann-Grad limit \cite{lang} is not appropriate in this regime.

  \subsection{Modifications of the microscopic model}
As explained above a modified macroscopic model (MA'), where the material coefficients $\herql$ depend on the relative velocity of the collision partners is more realistic. As the microscopic material coefficients $e_{rql} \in\nat_0$ are integer numbers they cannot depend on $(\xbo ,t)$, because they have to be Lipschitz continuous for our proof (see (C4)). So we define analogous to (\ref{adef}) and (\ref{jump}) a proper set of $m$ different, possible collision events
\bea
(\anrqs)^i :=\beta_{N,rq}^{k,i} \left(\int_0^t \chi_{M(N,r,s)} (k)\, a_{N,rq}^i (\xnk(s),s)\, ds\right) ,\ i=1,\ldots ,m,\nonumber
\eea
with rates
\bea
a_{N,rq}^i (\xnk(t),t)\! :=\!\min\Big\{ C_a\, ,\,\hat{a}_{rq}^i (\xnk(t),t)\Big((\snq(t) *W_N )(\xnk(t))-\delta_{r,q} W_N (\vec{0})/N \Big)\Big\}\nonumber
\eea
and corresponding outcomes $e_{rql}^i$. So the process with the most probable outcome for the relative velocity of the collision partners can be given the highest rate, whereas the others are small. With this microscopic model (MI') we introduce a dependence of the collision outcome on the relative velocity. To obtain compatibility with the macroscopic model (MA') certain conditions on $\hat{a}_{rq}^i$ and $e_{rql}^i$ have to be satisfied, and our proof of convergence applies with some minor changes.\\
Another modification of the microscopic model is to include pair interactions. Here we have to define a process for every pair of particles with $k\neq j$:
\bea
\anrqskj (t)&:=&\beta_{N,rq}^{kj}\left(\int_0^t \anrqkj (s)\, ds\right)\quad\,\mbox{for }k< j\nonumber\\
\anrqskj (t)&:=&\anqrsjk (t)\qquad\qquad\qquad\qquad\mbox{for }k>j,\nonumber
\eea
with rates
\bea
\anrqkj (t):=(1-\delta_{kj} )\chi_{M(N,r,t)} (k)\,\chi_{M(N,q,t)} (j)\,\oon\, W_N (\xnk(t)-\xnj(t))\,\harq (t).\nonumber
\eea
Our proof of convergence can be applied, but works only for spatially independent macroscopic collision rates $\harq(t)$. It assures convergence only up to a stopping time, as long as the summed rates $\sujrt\anrqkj (t)$ are bounded uniformly in $N$ by some predefined constant (cf.\ (M4)). For this reason the first modification seems to be more attractive, as it covers all realistic features explained in section V.A. Nevertheless pair interactions are more realistic descriptions of collision events and coagulation could be included in such a model.

  \subsection{Conclusion}
In this paper we specified a microscopic particle model (MI) from which we rigorously derived Smoluchowski's equation in the space dependent form (MA). Using the technique of the moderate limit developed in \cite{oel1}, we could estimate the critical fluctuation terms with martingale inequalities. This scaling assumption is a good approximiation of real systems in the regime of high temperatures and particle densities, which has been of interest in the study of interstellar gas clouds after shocks.\\
In the framework of the moderate limit, a further interesting question is the asymptotic behaviour of fluctuations for large system sizes and the formulation of a central limit theorem for this problem. One can as well try to eliminate some of the technical conditions for the proof, such as the restrictions to the scaling parameter in (C1) or the finite number of particle sizes. There is also hope to proof a convergence result for a microscopic model with pair interactions without the constraints mentioned in V.B, by using a suitable Sobolev-norm.\\
It would be certainly of most interest to derive Smoluchowski's equation in the hydrodynamic limit, but this task cannot be achieved with the methods used here. Nevertheless we could prove the validity of the spatially inhomogeneous equation in a regime, which is of great interest in astrophysics, and where previous derivations do not apply.

\end{document}